\documentclass[11pt]{amsart}

\theoremstyle{plain}
\newtheorem{theorem}                 {Theorem}      [section]
\newtheorem{proposition}  [theorem]  {Proposition}
\newtheorem{corollary}    [theorem]  {Corollary}
\newtheorem{lemma}        [theorem]  {Lemma}

\theoremstyle{definition}

\newtheorem{remark}       [theorem]  {Remark}
\newtheorem{definition}   [theorem]  {Definition}

\numberwithin{equation}{section}

\def \n{\mbox{${\mathbb N}$}}
\def \r{\mbox{${\mathbb R}$}}
\def \s{\mbox{${\mathbb S}$}}

\def \hor{\mbox{${\mathcal H}$}}
\def \ver{\mbox{${\mathcal V}$}}

\def \sq {\mbox{${\scriptstyle {\frac{r}{\sqrt{2}}}}$}}
\def \sqr {\mbox{${\scriptstyle {\frac{1}{\sqrt{2}}}}$}}

\def \tr{\mbox{${\tiny \r}$}}

\DeclareMathOperator{\dimension}{dim}
\DeclareMathOperator{\Span}{span}
\DeclareMathOperator{\trace}{trace}
\DeclareMathOperator{\grad}{grad}

\DeclareMathOperator{\Div}{div} 
 
\DeclareMathOperator{\ricci}{Ricci}
\DeclareMathOperator{\vol}{Vol}

\begin{document}

\title
{On the biharmonic and harmonic indices of the Hopf map}

\author{E. Loubeau}

\address{D\'epartement de Math\'ematiques\\
Laboratoire C.N.R.S. \ U.M.R. 6205\\
Universit\'e de Bretagne Occidentale\\
6, Avenue Victor Le Gorgeu\\
FRANCE}
\email{loubeau@univ-brest.fr}

\author{C. Oniciuc}

\address{Faculty of Mathematics\\
"Al.I.~Cuza" University of Iasi\\
Bd. Carol I, no. 11\\
700506 Iasi\\
ROMANIA}
\email{oniciucc@uaic.ro}



\thanks{
The authors are grateful to T.~Levasseur
for his help with representation theory. The second author
thanks the C.N.R.S. for a grant which made possible a three-month stay
at the Universit{\'e} de Bretagne Occidentale in Brest.}

\begin{abstract}
Biharmonic maps are the critical points of the bienergy functional and,
from this point of view, generalise harmonic maps. We consider the Hopf
map $\psi:\s^3\to \s^2$ and modify it into a nonharmonic biharmonic map
$\phi:\s^3\to \s^3$. We show $\phi$ to be unstable and estimate its
biharmonic index and nullity. Resolving the spectrum of the vertical
Laplacian associated to the Hopf map, we recover Urakawa's determination
of its harmonic index and nullity.
\end{abstract}

\maketitle

\section{Introduction}

In order to analyse the space of smooth maps between manifolds, Eells and
Sampson introduced in ~\cite{JEJHS} the family of functionals:
$$
E^k(\phi)=\frac{1}{2}\int_M\vert(d+d^{\star})^k\phi\vert^2 \ v_g,
\quad k\in \n^{\star}
$$
for maps $\phi:(M,g)\to (N,h)$.
\newline They quickly specialised their study to $k=1$ and called the
critical points of the {\em energy} $E^{1}$ {\em harmonic maps}, 
the corresponding
Euler-Lagrange equation being the vanishing of the tension field
$\tau(\phi)=\trace\nabla^{\phi}d\phi$. Here $\nabla^{\phi}$ denotes
the connection in the pull-back bundle $\phi^{-1}TN$.
\newline Recently, the case of the {\em bienergy} $E^2$ has been the subject
of some scrutiny. Its critical points are called {\em biharmonic maps}
and Jiang, in ~\cite{GYJ}, obtained the first and second variation formulas, 
showing, in particular, that a map is biharmonic if and only if
$\tau^2(\phi)=0$, where $\tau^2(\phi)=-J^{\phi}(\tau(\phi))$ and
$J^{\phi}$ is the Jacobi operator of the second variation
for the energy. As $J^{\phi}$ is
linear, any harmonic map is biharmonic, so we are interested in
nonharmonic biharmonic maps.
\newline Naturally, biharmonic submanifolds have been the first centre of
attention. On surfaces of revolution, Caddeo, Montaldo and Piu obtained
the parametric equation of all nongeodesic biharmonic curves
on the torus and Delaunay surfaces ~\cite{RCSMPP}, nonminimal biharmonic
submanifolds of $\s^3$ are classified in ~\cite{RCSMCO1}, while
constructions of such submanifolds on $\s^n$ ($n\geq 4$) are presented in
~\cite{RCSMCO2}.
\newline In the presence of nonconstant sectional curvature, a parametric
description of all nongeodesic biharmonic
curves of the Heisenberg group is given in ~\cite{RCCOPP}, whereas 
Inoguchi classified, in
~\cite{JI}, the biharmonic Legendre curves and Hopf cylinders on a
$3$-dimensional Sasakian space form and Sasahara gives in
~\cite{TS2} the explicit representation of nonminimal biharmonic
Legendre surfaces in a $5$-dimensional Sasakian space.
\newline While this results demonstrate the existence of nontrivial biharmonic
maps, though not in abundance, they carry little or no indication as to their
global behaviour. To this effect, some effort has been directed towards
the stability of biharmonic maps. As an example, the identity map
$\s^n\to \s^n$ and the totally geodesic inclusion
map $\s^m\to \s^n$ are biharmonic stable and their nullities
were computed in ~\cite{CO1}. In ~\cite{TS3}, the author proved that the
nonminimal biharmonic Legendre surfaces in a $5$-dimensional Sasakian
space are unstable.
\newline The index of the inclusion 
$\s^n(\sqr)\to \s^{n+1}$ was computed in ~\cite{ELCO}, and it is exactly 
one. Then, were investigated the indices of biharmonic maps in the unit
Euclidean sphere $\s^{n+1}$ obtained from minimal Riemannian immersions in 
$\s^n(\sqr)$. In particular, the authors showed that the index of the
nonminimal biharmonic map
$\s^m\big(\scriptstyle{\sqrt{\frac{m+1}{m}}}\big)\to \s^{m+p+1}$
derived from the generalised Veronese map, is at least $2m+3$.
\newline The present article pursues the same line of research
for submersions, by modifying a harmonic Riemannian submersion
$\psi:M\to \s^n(\sqr)$ into a nonharmonic
biharmonic subimmersion $\phi:M\to \s^{n+1}$ which we prove to be
unstable. Concentrating on the Hopf map
$\psi:\s^3(\sqrt{2})\to \s^2(\sqr)$, we first establish the
spectrum of the vertical Laplacian  of $\psi$ which allows us to
determine the harmonic index and nullity of $\psi$, hence recovering 
Urakawa's result (cf. ~\cite{HU}). We also give a geometric
characterization of $\ker J^{\psi}$, and obtain an
estimate of the biharmonic index and nullity of
$\phi:\s^3(\sqrt{2})\to \s^3$: the index is at least $11$ and the
nullity is greater than $8$.

Throughout the paper, manifolds, metrics, maps are assumed to be smooth, 
and $(M,g)$ denotes a connected Riemannian manifold without boundary. We 
denote by $\nabla$ the Levi-Civita connection of $(M,g)$, and for
the curvature tensor field $R$ we use the sign convention 
$R(X,Y)=[\nabla_X,\nabla_Y]-\nabla_{[X,Y]}$, while the Laplacian on 
sections of $\phi^{-1}TN$ is 
$\Delta^{\phi}=-\trace (\nabla^{\phi})^2$.  

Let
$$
\s^n(\sq)=\s^n(\sq)\times\{\sq\}=
\{p=(x^1,\ldots,x^{n+1},\sq)\big{\vert}
(x^1)^2+\cdots+(x^{n+1})^2=\scriptstyle{\frac{r^2}{2}}\}
$$
be a hypersphere of $\s^{n+1}(r)$, as proved in
~\cite{RCSMCO1}, $\s^n(\sq)$ is a nonminimal biharmonic
submanifold of $\s^{n+1}(r)$. Indeed, let
$\eta=\frac{1}{r}(x^1,\ldots,x^{n+1},-\sq)$ be a unit section of
the normal bundle of $\s^n(\sq)$ in $\s^{n+1}(r)$. Then, the
second fundamental form of $\s^n(\sq)$ is $B(X,Y)=\nabla d{\bf
i}(X,Y)=-\frac{1}{r}\langle X,Y\rangle\eta$ and the tension field
of the inclusion map ${\bf i}:\s^n(\sq)\to \s^{n+1}(r)$ is
$\tau({\bf i})=-\frac{n}{r}\eta$. Besides, by direct computation,
$\tau^2({\bf i})=0$.

\begin{theorem}
Let $M$ be a compact manifold and $\psi:M\to \s^n(\sq)$ a
nonconstant map. The map $\phi={\bf i}\circ \psi:M\to \s^{n+1}(r)$
is nonharmonic biharmonic if and only if $\psi$ is harmonic and
$e(\psi)$ is constant.
\end{theorem}

\begin{proof}
The composition law gives
$$
\tau(\phi)=\tau(\psi)-\frac{2}{r}e(\psi)\eta.
$$
By straightforward calculations:
\begin{eqnarray*}
\Delta^{\phi}\tau(\phi)&=&\Delta^{\psi}\tau(\psi)
+\frac{2}{r^2}e(\psi)\tau(\psi)+\frac{1}{r^2}d\phi(\theta^{\sharp}+
4\grad e(\psi)) \\
&&-\frac{1}{r}\big( \frac{4}{r^2}(e(\psi))^2+2\Delta e(\psi)
-2\Div\theta^{\sharp}+\vert\tau(\psi)\vert^2 \big)\eta,
\end{eqnarray*}
and
$$
\trace R^{\tiny{\s^{n+1}(r)}}(d\phi\cdot,\tau(\phi))d\phi\cdot=
\frac{1}{r^2}d\phi(\theta^{\sharp})-\frac{2}{r^2}e(\psi)\tau(\psi)+
\frac{4}{r^3}(e(\psi))^2\eta,
$$
where $\theta$ is a $1$-form on $M$ given by
$\theta(X)=\langle d\psi(X),\tau(\psi)\rangle$, and
$\theta^{\sharp}\in C(TM)$ is defined by
$\langle \theta^{\sharp},X\rangle=\theta(X)$, $X\in C(TM)$.

\noindent
Replacing in the biharmonic equation we obtain:
$$
\tau^2(\phi)=-\Delta^{\psi}\tau(\psi)-\frac{2}{r^2}d\phi(\theta^{\sharp}
+2\grad e(\psi))
+\frac{1}{r}\big( 2\Delta e(\psi)-2\Div\theta^{\sharp}
+\vert\tau(\psi)\vert^2 \big)\eta.
$$

Assume that $\phi$ is biharmonic. As the $\eta$-part of $\tau^2(\phi)$
vanishes:
$$
2\Delta e(\psi)-2\Div\theta^{\sharp}+\vert\tau(\psi)\vert^2=0.
$$
By Stokes, we deduce that $\tau(\psi)=0$ which
implies $\Delta e(\psi)=0$, i.e.
$e(\psi)$ is constant.

The converse is immediate.
\end{proof}

\begin{corollary}[\cite{RCSMCO2,CO2}]
If $\psi$ is a harmonic Riemannian immersion or submersion, then $\phi$ is
nonharmonic biharmonic.
\end{corollary}

We close this section with the second variation formula for
biharmonic maps in $\s^n(r)$.

\begin{theorem}
[\cite{CO1}] Let $\phi:(M,g)\to \s^n(r)$ be a biharmonic map. The
Hessian of the bienergy $E^2$ at $\phi$ is:
$$
H(E^2)_{\phi}(V,W)=\int_M\langle I^{\phi}(V),W\rangle \ v_g,
$$
with
\begin{eqnarray}
\label{eq:J1}
r^2I^{\phi}(V)&=&r^2\Delta^{\phi}(\Delta^{\phi}
V)+\Delta^{\phi}\big( \trace\langle V,d\phi\cdot\rangle d\phi\cdot-
\vert d\phi\vert^2V \big)
+2\langle d\tau(\phi),d\phi\rangle V
\nonumber\\
&&+\vert\tau(\phi)\vert^2V
-2\trace\langle V,d\tau(\phi)\cdot\rangle d\phi\cdot
-2\trace\langle\tau(\phi),dV\cdot\rangle d\phi\cdot
\nonumber\\
&&-\langle\tau(\phi),V\rangle\tau(\phi)+\trace\langle
d\phi\cdot,\Delta^{\phi} V\rangle d\phi\cdot
+\frac{1}{r^2}\trace\big\langle d\phi\cdot,\trace\langle V,d\phi\cdot\rangle
d\phi\cdot\big\rangle d\phi\cdot
\\
&&-\frac{2}{r^2}\vert d\phi\vert^2\trace\langle d\phi\cdot,
V\rangle d\phi\cdot
+2\langle dV,d\phi\rangle\tau(\phi)-\vert
d\phi\vert^2\Delta^{\phi} V+ \frac{1}{r^2}\vert d\phi\vert^4 V,
\nonumber
\end{eqnarray}
and $V,W\in C(\phi^{-1}T\s^n(r))$.
\end{theorem}

\section{Some general properties of the hopf map $\psi$}

Let $\psi:(M,g)\to (N,h)$ be a submersion. For $p\in M$, call
$\ver_p=\ker d\psi_p$ the {\em vertical space} and $\hor_p$, its
orthogonal complement, the {\em horizontal space}. Any $X\in T_pM$ splits into
$X^{\tiny{\ver}}+X^{\tiny{\hor}}$, where $X^{\tiny{\ver}}\in
\ver_p$ and $X^{\tiny{\hor}}\in \hor_p$, and a vector field $X$ is
{\em basic} if horizontal, i.e.
$X=X^{\tiny{\hor}}$, and $\psi$-related to
$Y\in C(TN)$, i.e. $d\psi_p(X(p))=Y(\psi(p))$, $\forall p\in M$.

The map $\psi$ is a {\em Riemannian submersion} if:
$$
h(d\psi(X),d\psi(Y))=g(X,Y), \quad \forall X,Y\in \hor_p.
$$

Let $\psi:\s^3(\sqrt{2})\to \s^2(\sqr)$ be the Hopf map given by:
$$
\psi(z,w)=\frac{1}{2\sqrt{2}}(2z\overline{w},\vert z\vert^2-\vert
w\vert^2),
$$
or
$$
\psi(x^1,x^2,x^3,x^4)=\frac{1}{2\sqrt{2}}(2x^1x^3+2x^2x^4,
2x^2x^3-2x^1x^4,(x^1)^2+(x^2)^2-(x^3)^2-(x^4)^2),
$$
where $z=x^1+{\mathrm i}x^2$ and $w=x^3+{\mathrm i}x^4$. The map
$\psi$ is a harmonic Riemannian submersion and its fibres are
totally geodesic submanifolds of $\s^3(\sqrt{2})$. More precisely,
$\psi^{-1}\big(\psi(z,w)\big)$ is the great circle passing through
$(z,w)$ and ${\mathrm i}(z,w)=({\mathrm i}z,{\mathrm i}w)$. Remark
that $\psi^1$, $\psi^2$ and $\psi^3$ are harmonic homogeneous
polynomials of degree $2$ on $\r^4$.

Let $x\in \s^3(\sqrt{2})$, $y\in \s^2(\sqr)$ and define the vector
fields:
$$
X_1(x)=\frac{1}{\sqrt{2}}(-x^2,x^1,-x^4,x^3), \quad
X_2(x)=\frac{1}{\sqrt{2}}(-x^3,x^4,x^1,-x^2),
$$
$$
X_3(x)=\frac{1}{\sqrt{2}}(-x^4,-x^3,x^2,x^1), \quad
X_4(x)=\frac{1}{\sqrt{2}}(-x^2,x^1,x^4,-x^3),
$$
$$
X_5(x)=\frac{1}{\sqrt{2}}(-x^3,-x^4,x^1,x^2), \quad
X_6(x)=\frac{1}{\sqrt{2}}(-x^4,x^3,-x^2,x^1),
$$
$$
Y_1(y)=\sqrt{2}(-y^2,y^1,0), \quad Y_2(y)=\sqrt{2}(y^3,0,-y^1), \quad
Y_3(y)=\sqrt{2}(0,-y^3,y^2).
$$
Then:
\begin{itemize}
\item
$\{X_1,\ldots ,X_6\}$ is a basis of Killing vector fields on $\s^3(\sqrt{2})$,
\item
$\{Y_1,Y_2,Y_3\}$ is a basis of Killing vector fields on
$\s^2(\sqr)$,
\item
$X_4$ is $\psi$-related to $Y_1$, $X_5$ is
$\psi$-related to $Y_2$, and $X_6$ is $\psi$-related to $Y_3$,
\item
$\{X_1,X_2,X_3\}$ is a global orthonormal frame on $\s^3(\sqrt{2})$,
\item
$X_1=X_1^{\tiny{\ver}}$ and $\ver\s^3(\sqrt{2})=\Span\{X_1\}$,
\item
$X_2=X_2^{\tiny{\hor}}$, $X_3=X_3^{\tiny{\hor}}$ and
$\hor\s^3(\sqrt{2})=\Span\{X_2,X_3\}$,
\item
$d\psi_x(X_2)$ and $d\psi_x(X_3)$ form an orthonormal basis of
$T_{\psi(x)}\s^2(\sqr)$.
\end{itemize}
A straightforward computation yields:

\begin{proposition}
\label{eq:connection}
The Levi-Civita connection on $\s^3(\sqrt{2})$ is given by:
$$
\left\{
\begin{array}{l}
\nabla_{X_1}X_1=0, \quad \nabla_{X_1}{X_2}=-\frac{1}{\sqrt{2}}X_3,
\quad \nabla_{X_1}X_3=\frac{1}{\sqrt{2}}X_2,
\\ \mbox{} \\
\nabla_{X_2}X_1=\frac{1}{\sqrt{2}}X_3, \quad \nabla_{X_2}X_2=0,
\quad \nabla_{X_2}X_3=-\frac{1}{\sqrt{2}}X_1,
\\ \mbox{} \\
\nabla_{X_3}X_1=-\frac{1}{\sqrt{2}}X_2, \quad
\nabla_{X_3}X_2=\frac{1}{\sqrt{2}}X_1, \quad \nabla_{X_3}X_3=0,
\end{array}
\right.
$$
and
$$
[X_1,X_2]=-\sqrt{2}X_3, \quad
[X_2,X_3]=-\sqrt{2}X_1, \quad [X_3,X_1]=-\sqrt{2}X_2.
$$
\end{proposition}

Recall two general properties of Killing vector
fields which will be used throughout the paper:

\noindent $i)$ Let $X$ be a Killing
vector field on a compact Riemannian manifold $(M,g)$. Denote
by $\{\varphi_t\}_{t\in \tr}$ its flow, and consider $f\in
C^{\infty}(M)$, then, as $\varphi_t$ is an isometry:
$$
\int_Mf \ v_g=\int_M(f\circ\varphi_t) \ v_g, \quad \forall t\in \r,
$$
and deriving by $t$ produces:
$$
0=\frac{d}{dt}\Big\vert_{t=0}\int_M(f\circ\varphi_t) \ v_g=\int_M(Xf) \ v_g.
$$
Alternatively, if $X$ is Killing, then $Xf=\Div (fX)$ and, by Stokes,
$\int_M(Xf) \ v_g=0$. In particular:
$$
\int_M(Xf_1)f_2 \ v_g=-\int_Mf_1(Xf_2) \ v_g
\quad \hbox{and} \quad
\int_M(Xf)f \ v_g=0.
$$

\noindent $ii)$ Having an isometry in $\varphi_t$, implies
$\Delta(f\circ\varphi_t)=(\Delta f)\circ\varphi_t$, for any $t$,
and, again, deriving by $t$, means that:
$$
\Delta(Xf)=X(\Delta f),
$$
so, if $\Delta f=\lambda f$, $\Delta(Xf)=\lambda(Xf)$, i.e. $X$
preserves the eigenspaces of the Laplacian.

More specific to the Hopf map is:

\begin{lemma}
\label{eq:secondlemma}
$$
\dimension\{f\in C^{\infty}(\s^3(\sqrt{2}))\big\vert X_1f=0, \
\Delta f=\lambda_{2k}f\}= \dimension\{\tilde{f}\in
C^{\infty}(\s^2(\sqr))\big\vert \Delta\tilde{f}=\mu_k\tilde{f}\},
$$
where $\lambda_k=\frac{k(k+2)}{2}$ and $\mu_l=2l(l+1)$, $k,l\in
\n$, are the eigenvalues of the Laplacians of $\s^3(\sqrt{2})$ and
$\s^2(\sqr)$.
\end{lemma}

\begin{proof}
As $\psi$ is a Riemannian submersion and its fibres are totally
geodesic (hence minimal)
submanifolds of $\s^3(\sqrt{2})$ (see ~\cite{MBPGEM}):
$$
\Delta(\tilde{f}\circ\psi)=(\Delta\tilde{f})\circ\psi, \quad
\forall \tilde{f}\in C^{\infty}(\s^2(\sqr)).
$$
Therefore
$$
T:\{\tilde{f}\in C^{\infty}(\s^2(\sqr))\big\vert
\Delta\tilde{f}=\mu_k\tilde{f}\} \to \{f\in
C^{\infty}(\s^3(\sqrt{2}))\big\vert X_1f=0, \ \Delta
f=\lambda_{2k}f\}
$$
defined by
$
T(\tilde{f})=\tilde{f}\circ\psi
$
is an isomorphism and the lemma follows.
\end{proof}

Like any Riemannian submersion, $\psi:\s^3(\sqrt{2})\to
\s^2(\sqr)$ satisfies the coarea formula (cf. ~\cite{F}):
$$
\int_{\s^3(\sqrt{2})}(\tilde{f}\circ\psi)
 \ v_g=2\sqrt{2}\pi\int_{\s^2(\sqr)}\tilde{f} \ v_g, \quad
\forall \tilde{f}\in C^{\infty}(\s^2(\sqr)).
$$

\section{The vertical laplacian}

To a Riemannian submersion $\psi:(M,g)\to (N,h)$ one associates a
second-order operator by restricting functions to a fibre and, viewing it
as a submanifold, consider its induced Laplacian.

\begin{definition}[\cite{LBBJPB}]
The {\em vertical Laplacian} $\Delta^V$ is the differential
operator defined on $(M,g)$ by:
$$
(\Delta^Vf)(p)=\big(\Delta^{F_p}(f_{\vert F_p})\big)(p), \quad
f\in C^{\infty}(M),
$$
where $\Delta^{F_p}$ is the Laplace operator of the fibre
$F_p=\psi^{-1}(\psi(p))$.
\end{definition}

\noindent
The fibres being isometric, the spectrum of $\Delta^V$ is discrete.

Berard Bergery and Bourguignon proved in ~\cite{LBBJPB}
that, when the fibres are totally geodesic, the Laplacian $\Delta$ on $M$
and the vertical Laplacian commute, implying:

\begin{theorem}[\cite{LBBJPB}]
\label{eq:Hilbertbasis}
The Hilbert space $L^2(M)$ admits a Hilbert
basis consisting of simultaneous eigenfunctions for $\Delta$ and
$\Delta^V$.
\end{theorem}

Denote by
$\{c_i\}_{i\in \tiny{\n}}$ the eigenvalues of $\Delta^V$
(shared with the Laplacian
$\Delta^F$ of a generic fibre $F$),
then $c_0=0$, $c_i>0$ for any $i>0$, and
$c_i\longrightarrow\infty$ when $i\longrightarrow\infty$. Note
that the multiplicities are not necessarily finite.

For the Hopf map, the corresponding vertical Laplacian is
$\Delta^Vf=-X_1X_1f$. As $X_1$, and by consequence $\Delta^V$, preserves the
eigenspaces of $\Delta$, to establish the spectrum of the
vertical Laplacian we restrict $\Delta^V$ to the eigenspaces of $\Delta$,
and consider the following problem:
\begin{equation}
\label{eq:vertical1}
\left\{
\begin{array}{l}
\Delta f=\lambda_kf
\\ \mbox{} \\
\Delta^Vf=c_lf.
\end{array}
\right.
\end{equation}
We prove that $c_{l}=\frac{(k-2l)^2}{2}$, $l\in \left\{0,
\dots,k\right\}$ with multiplicity $2(k+1)$,
with the exception of
$c_{k/2}=0$, when $k$ is even, which has multiplicity $k+1$:

Denote by ${\mathbf H}^{k}$ the set of harmonic homogeneous
polynomials of degree $k$ on ${\mathbb C\,}^{2}$ and recall that
$\dimension_{{\mathbb C\,}}{{\mathbf H}^{k}} = (k+1)^{2}$. We work
in the algebra of differential operators on polynomials
$R={\mathbb C\,}[z,\bar{z},w,\bar{w}]$ and introduce the linear
operators:
$$
\Delta = -4 (\partial_{z}\partial_{\bar{z}} +
\partial_{w}\partial_{\bar{w}}) , \quad {\mathsf h} =
\bar{z}\partial_{\bar{z}} - z\partial_{z} +
\bar{w}\partial_{\bar{w}} - w\partial_{w} .
$$
Note that $\Delta$ is the usual Laplacian on ${\mathbb C\,}^{2}$
and:
$$
[\Delta , {\mathsf h}] = 0,
$$
therefore ${\mathbf H}^{k}$ is stable under the action of
${\mathsf h}$.
\newline Since ${\mathsf h} = {\mathrm i} \sqrt{2} X_{1}$, we
determine the spectrum of ${\mathsf h}$ restricted to ${\mathbf H}^{k}$
to solve ~\eqref{eq:vertical1}.
The strategy is to include ${\mathsf h}$ in a Lie algebra isomorphic
to $\mathfrak{sl}(2,{\mathbb C\,})$ acting on ${\mathbf H}^{k}$ and
use representation theory to deduce its spectrum.
\newline Let:
$$
{\mathsf e} = {\mathrm i}(\bar{w}\partial_{z} - \bar{z}\partial_{w}) ,
\quad {\mathsf f} = {\mathrm i}(w\partial_{\bar{z}} -
z\partial_{\bar{w}}).
$$
One can easily check that:
$$
[{\mathsf e},{\mathsf f}]={\mathsf h}  , \quad
[{\mathsf h},{\mathsf e}]=2{\mathsf e}  , \quad [{\mathsf
h},{\mathsf f}]=-2{\mathsf f}  ,
$$
making ${\mathfrak g} =
{\mathbb C\,}{\mathsf e} + {\mathbb C\,}{\mathsf f} + {\mathbb
C\,}{\mathsf h}$ a Lie subalgebra of ${\mathrm End}_{\mathbb
C\,}R$ isomorphic to $\mathfrak{sl}(2,{\mathbb C\,})$. Moreover
$$
[\Delta , {\mathsf e}] = [\Delta , {\mathsf f}] = [\Delta ,
{\mathsf h}] = 0,
$$
so ${\mathbf H}^{k}$ is stable under the action of ${\mathfrak g}$
and we can  decompose it as a ${\mathfrak g}$-module.
\newline Recall (cf.~\cite{WFJH}) that there exists, up to
isomorphisms, a unique
${\mathfrak g}$-module of dimension $d+1$:
$$
E(d) = \bigoplus_{l=0}^{d} {\mathbb C\,}{\mathsf f}^{l}(v),
$$
where $v\in E(d)$ is determined, up to a non-zero scalar, by:
$$
{\mathsf e}(v) = 0, \quad {\mathsf h}(v) = d v .
$$
This space $E(d)$ is called the simple ${\mathfrak g}$-module of
highest weight $d$ and $v$ a highest weight vector. Besides, each
${\mathbb C\,}{\mathsf f}^{l}(v)$ is the eigenspace of ${\mathsf
h}$ associated to the eigenvalue $d-2l$ ($l\in \left\{0,
\dots,d\right\}$), i.e.
$$
{\mathsf h} ({\mathsf f}^{l}(v)) = (d-2l) {\mathsf f}^{l}(v) .
$$
For $n \in \left\{0, \dots,k\right\}$, let $f_{n} = {\bar{z}}^{n}
{\bar{w}}^{k-n}$, then:
$$
\Delta(f_{n}) = 0, \quad {\mathsf e}(f_{n}) = 0,
\quad {\mathsf h}(f_{n}) = k f_{n} ,
$$
so $f_{n} \in {\mathbf H}^{k}$ spans a simple ${\mathfrak
g}$-submodule of highest weight $k$, denoted by $V(n)$, that is $V(n)
\simeq E(k)$ as ${\mathfrak g}$-modules. We deduce that $W=
\sum_{n=0}^{k} V(n)$ is a ${\mathfrak g}$-submodule of ${\mathbf
H}^{k}$. To show that this sum is direct, and infer from
$\dimension{V(n)} = k+1$ that ${\mathbf H}^{k} =
\bigoplus_{n=0}^{k} V(n)$, we consider the operator:
$$
\Lambda = (z\partial_{z} - \bar{z}\partial_{\bar{z}})
- (w\partial_{w} - \bar{w}\partial_{\bar{w}}) .
$$
Observe that:
$$
[\Lambda , {\mathsf f}] = 0 , \quad \Lambda (f_{n}) = (k-2n) f_{n},
\quad \Lambda ({\mathsf f}^{l}(f_{n}))
= {\mathsf f}^{l} (\Lambda(f_{n}))= (k-2n) {\mathsf f}^{l}(f_{n}).
$$
Therefore $V(n)$ is the eigenspace of $\Lambda$ associated to the
eigenvalue $(k-2n)$, and ${\mathbf H}^{k}$ is the direct sum of
the $V(n)$'s.
\newline This shows that the eigenvalues of
${\mathsf h}$ are $k-2l$, $l \in \{0,\dots,k\}$ with multiplicity $k+1$
and eigenspaces $\bigoplus_{n=0}^{k} {\mathbb C\,}{\mathsf f}^{l}(f_{n})$.
\newline In conclusion, the spectrum of the vertical Laplacian
$\Delta^{V} = - X_{1} X_{1} = \frac{1}{2} {\mathsf h}^{2}$ is
$\left\{ \frac{(k-2l)^2}{2} : l \in \{1,\dots,k\} \right\}$
with multiplicity $2(k+1)$, with, when $k$ is even, $0$ of multiplicity
$(k+1)$.

\section{The harmonic index and nullity of $\psi$}

The second variation of a harmonic map $\psi:(M,g)\to (N,h)$ is given by
the Jacobi operator $J^{\psi}$ (cf. ~\cite{EM,RTS}):
$$
H(E^1)_{\psi}(V,W)=\int_M\langle J^{\psi}(V),W\rangle \ v_g,
$$
where $V,W\in C(\psi^{-1}TN)$ and $J^{\psi}=\Delta^{\psi}+ \trace
R^{\tiny{N}}(d\psi, \ )d\psi$. In the case of the Hopf map, the
Jacobi operator $J^{\psi}:C(\psi^{-1}T\s^2(\sqr))\to
C(\psi^{-1}T\s^2(\sqr))$ is:
\begin{eqnarray*}
J^{\psi}(V)&=&\Delta^{\psi}V+\trace
R^{\tiny{\s^2(\sqr)}}(d\psi,V)d\psi \\
&=&\Delta^{\psi}V-\ricci^{\tiny{\s^2(\sqr)}}(V) \\
&=&\Delta^{\psi}V-2V.
\end{eqnarray*}

As $\{X_2(x),X_3(x)\}$ is a basis of $\hor_x$, $\forall x\in
\s^3(\sqrt{2})$, any section $V$ can be written
$f_2d\psi(X_2)+f_3d\psi(X_3)$, where $f_2,f_3\in
C^{\infty}(\s^3(\sqrt{2}))$, hence
$$
C(\psi^{-1}T\s^2(\sqr))=\{f_2d\psi(X_2)\vert f_2\in
C^{\infty}(\s^3(\sqrt{2}))\} \oplus \{f_3d\psi(X_3)\vert f_3\in
C^{\infty}(\s^3(\sqrt{2}))\}.
$$
This sum is orthogonal with respect to the usual scalar product on
$C(\psi^{-1}T\s^2(\sqr))$.

A simple computation results in:

\begin{lemma}
\label{eq:Jlemma1}
If $X\in C(T\s^3(\sqrt{2}))$,
then:
$$
\Delta^{\psi}(d\psi(X))=d\psi(X-\trace\nabla^2X)-2\trace\nabla^{\psi}
d\psi( \ ,\nabla X).
$$
\end{lemma}

Using Proposition ~\ref{eq:connection}, Lemma ~\ref{eq:Jlemma1},
the features of Killing
vector fields and the properties of the second fundamental form of a
Riemannian submersion, we have:

\begin{lemma}
\label{eq:Jlemma2} The action of $J^{\psi}$ on
$C(\psi^{-1}T\s^2(\sqr))$ is specified by:
$$
J^{\psi}(f_2d\psi(X_2))=(\Delta
f_2)d\psi(X_2)+2\sqrt{2}(X_1f_2)d\psi(X_3)
$$
and
$$
J^{\psi}(f_3d\psi(X_3))=(\Delta
f_3)d\psi(X_3)-2\sqrt{2}(X_1f_3)d\psi(X_2)
$$
\end{lemma}

By Lemma ~\ref{eq:Jlemma2}, if $\Delta f_2=\lambda_kf_2$ and
$\Delta f_3=\lambda_kf_3$, then:
$$
J^{\psi}(f_2d\psi(X_2))=\lambda_k
f_2d\psi(X_2)+2\sqrt{2}(X_1f_2)d\psi(X_3),
$$
$$
J^{\psi}(f_3d\psi(X_3))=\lambda_k
f_3d\psi(X_3)-2\sqrt{2}(X_1f_3)d\psi(X_2).
$$
Moreover, since $X_1$ preserves the eigenspaces of the Laplacian,
$J^{\psi}$ leaves invariant the subspace
$$
S^{\psi}_{\lambda_k}=\{f_2d\psi(X_2)\vert \Delta
f_2=\lambda_kf_2\}\oplus \{f_3d\psi(X_3)\vert \Delta
f_3=\lambda_kf_3\},
$$
for any $k\in \n$. If $k_1\neq k_2$,
$S^{\psi}_{\lambda_{k_1}}$ and $S^{\psi}_{\lambda_{k_2}}$ are orthogonal one
to the other.

\begin{proposition}
The eigenvalues of $J^{\psi}$ restricted to $S^{\psi}_{\lambda_k}$
are $\lambda_k\pm 2\sqrt{2c}$, where $c$ is an eigenvalue of
$\Delta^V$ restricted to $\{f\in C^{\infty}(\s^3(\sqrt{2}))\vert
\Delta f=\lambda_k f\}$, and the eigensections are:
\begin{itemize}
\item if $c=0$ and $X_1f_2=X_1f_3=0$, then
$$
J^{\psi}(f_2d\psi(X_2)+f_3d\psi(X_3))=\lambda_k\big(f_2d\psi(X_2)
+f_3d\psi(X_3)\big),
$$
\item if $c>0$ and $\Delta^Vf=cf$, then
$$
J^{\psi}(fd\psi(X_2)+\frac{1}{\sqrt{c}}(X_1f)d\psi(X_3))=(\lambda_k
+2\sqrt{2c})
\big(fd\psi(X_2)+\frac{1}{\sqrt{c}}(X_1f)d\psi(X_3)\big)
$$
and
$$
J^{\psi}(fd\psi(X_2)-\frac{1}{\sqrt{c}}(X_1f)d\psi(X_3))=(\lambda_k
-2\sqrt{2c})
\big(fd\psi(X_2)-\frac{1}{\sqrt{c}}(X_1f)d\psi(X_3)\big).
$$
\end{itemize}
The multiplicity of each eigenvalue of $J^{\psi}$ restricted to
$S^{\psi}_{\lambda_k}$ is $2(k+1)$.
\end{proposition}

The spectrum of the vertical Laplacian (consult Section $3$) means that:

If $k=0$, then $\lambda_0=0$, $c=0$ and $\lambda_k\pm 2\sqrt{2c}=0$.

If $k=1$, then $\lambda_1=\frac{3}{2}$, $c=\frac{1}{2}$ (with multiplicity
$4$) thus $\lambda_k\pm 2\sqrt{2c}\in
\{-\frac{1}{2},\frac{7}{2}\}$.

If $k=2$, then $\lambda_2=4$, $c=0$ or $c=2$ and
$\lambda_k\pm 2\sqrt{2c}\in \{0,4,8\}$. The multiplicity of $c=0$
is $3$ and the multiplicity of $c=2$ is $6$.

Notice that $c=\frac{(k-2l)^2}{2}\leq\frac{k^2}{2}$ so, if
$k>2$, $\lambda_k>2\sqrt{2c}$, for any $c$, and the eigenvalues of
$J^{\psi}$ are positive.

Homothetic transformations do not change the spectrum of the operator
$J^{\psi}$  thus, from the above
analysis, we recover Urakawa's result:

\begin{theorem}[\cite{HU}]
\label{eq:Urakawa}
The Hopf map $\psi:\s^3\to \s^2$ is harmonic of index $4$ and nullity $8$.
\end{theorem}

An alternative proof of Theorem ~\ref{eq:Urakawa} is via a direct method:

\noindent As $J^{\psi}$ preserves $S^{\psi}_{\lambda_k}$,
we find an $L^2$-orthonormal basis of $S^{\psi}_{\lambda_k}$, $k=0$, $1$ or
$2$, and compute the matrix of $J^{\psi}$ in this basis. Then,
for $k>2$, we show that $J^{\psi}$
restricted to $S^{\psi}_{\lambda_k}$ is positive definite. The
index and nullity of the Hopf map will be given by the matrices
of the operator $J^{\psi}$ restricted to $S^{\psi}_{\lambda_k}$:

$i) \ k=0.$ In this case $f_2$ and $f_3$ are constants, and an
$L^2$-orthonormal basis of $S^{\psi}_{\lambda_0}$ is
$$
B^{\psi}_{\lambda_0}=\big\{\frac{1}{c}d\psi(X_2),\frac{1}{c}d\psi(X_3)\big\},
$$
where $c^2=\vol(\s^3(\sqrt{2}))$. We have
$J^{\psi}(d\psi(X_2))=J^{\psi}(d\psi(X_3))=0$, so
$J^{\psi}$, in the basis $B^{\psi}_{\lambda_0}$, is the zero matrix.

$ii) \ k=1.$ A basis of the first eigenfunctions of the Laplacian on
$\s^3(\sqrt{2})$ is $\{f^1_i(x)=\frac{x^i}{a}\}_{i=1}^4$, where
$2a^2=\vol\big(\s^3(\sqrt{2})\big)$, making
$$
B^{\psi}_{\lambda_1}=\{f^1_id\psi(X_2)\}_{i=1}^4\cup
\{f^1_jd\psi(X_3)\}_{j=1}^4
$$
an $L^2$-orthonormal basis for $S^{\psi}_{\lambda_1}$.
The matrix of $J^{\psi}$ restricted to
$S^{\psi}_{\lambda_1}$ is conditioned by the action of $X_1$, $X_2$ and
$X_3$ on $B^{\psi}_{\lambda_1}$:
$$
X_1f^1_1=-\frac{f^1_2}{\sqrt{2}}, \
X_1f^1_2=\frac{f^1_1}{\sqrt{2}}, \
X_1f^1_3=-\frac{f^1_4}{\sqrt{2}}, \
X_1f^1_4=\frac{f^1_3}{\sqrt{2}},
$$
$$
X_2f^1_1=-\frac{f^1_3}{\sqrt{2}}, \
X_2f^1_2=\frac{f^1_4}{\sqrt{2}}, \
X_2f^1_3=\frac{f^1_1}{\sqrt{2}}, \
X_2f^1_4=-\frac{f^1_2}{\sqrt{2}},
$$
$$
X_3f^1_1=-\frac{f^1_4}{\sqrt{2}}, \
X_3f^1_2=-\frac{f^1_3}{\sqrt{2}}, \
X_3f^1_3=\frac{f^1_2}{\sqrt{2}}, \
X_3f^1_4=\frac{f^1_1}{\sqrt{2}}.
$$
Calculating $J^{\psi}(f^1_id\psi(X_2))$ and
$J^{\psi}(f^1_jd\psi(X_3))$ yields the matrix:
$$
\left(
\begin{array}{cccccccccccc}
\frac{3}{2} & 0 & 0 & 0 & 0 & -2 & 0 & 0 \\
0 & \frac{3}{2} & 0 & 0 & 2 & 0 & 0 & 0 \\
0 & 0 & \frac{3}{2} & 0 & 0 & 0 & 0 & -2 \\
0 & 0 & 0 & \frac{3}{2} & 0 & 0 & 2 & 0 \\
0 & 2 & 0 & 0 & \frac{3}{2} & 0 & 0 & 0 \\
-2 & 0 & 0 & 0 & 0 & \frac{3}{2} & 0 & 0 \\
0 & 0 & 0 & 2 & 0 & 0 & \frac{3}{2} & 0 \\
0 & 0 & -2 & 0 & 0 & 0 & 0 & \frac{3}{2}
\end{array}
\right).
$$
Its eigenvalues are $-\frac{1}{2}$ and $\frac{7}{2}$, both of
multiplicity $4$.

$iii) \ k=2.$ Set:
$$
f^2_1=\frac{1}{b}(x^1x^2+x^3x^4), \quad
f^2_2=\frac{1}{b}(x^1x^2-x^3x^4), \quad
f^2_3=\frac{1}{b}(x^1x^3+x^2x^4),
$$
$$
f^2_4=\frac{1}{b}(x^1x^3-x^2x^4), \quad
f^2_5=\frac{1}{b}(x^1x^4+x^2x^3), \quad
f^2_6=\frac{1}{b}(x^1x^4-x^2x^3),
$$
$$
f^2_7=\frac{1}{2b}\big((x^1)^2+(x^2)^2-(x^3)^2-(x^4)^2\big), \quad
f^2_8=\frac{1}{2b}\big((x^1)^2-(x^2)^2+(x^3)^2-(x^4)^2\big),
$$
$$
f^2_9=\frac{1}{2b}\big((x^1)^2-(x^2)^2-(x^3)^2+(x^4)^2\big),
$$
where $3b^2=\vol\big(\s^3(\sqrt{2})\big)$, then $\{f^2_{i}\}_{i=1}^9$ is an
$L^2$-orthonormal basis of
$\{f\in C^{\infty}(\s^3(\sqrt{2}))\vert\Delta f=\lambda_2 f\}$,
and
$$
B^{\psi}_{\lambda_2}=\{f^2_id\psi(X_2)\}_{i=1}^9\cup
\{f^2_jd\psi(X_3)\}_{j=1}^9
$$
an $L^2$-orthonormal basis for $S^{\psi}_{\lambda_2}$. The
vector fields $X_1$, $X_2$, $X_3$ act upon $B^{\psi}_{\lambda_2}$ by:
\begin{eqnarray*}
\left\{
\begin{array}{l}
X_1f^2_1=\sqrt{2}f^2_8, \
X_1f^2_2=\sqrt{2}f^2_9, \ X_1f^2_3=0, \
X_1f^2_4=-\sqrt{2}f^2_5, X_1f^2_5=\sqrt{2}f^2_4,
\\ \mbox{} \\
X_1f^2_6=0, \
X_1f^2_7=0, \ X_1f^2_8=-\sqrt{2}f^2_1, \ X_1f^2_9=-\sqrt{2}f^2_2,
\end{array}
\right.
\end{eqnarray*}

\begin{eqnarray*}
\left\{
\begin{array}{l}
X_2f^2_1=\sqrt{2}f^2_6, \ X_2f^2_2=0, \ X_2f^2_3=\sqrt{2}f^2_9, \
X_2f^2_4=\sqrt{2}f^2_7, \ X_2f^2_5=0,
\\ \mbox{} \\
X_2f^2_6=-\sqrt{2}f^2_1, \ X_2f^2_7=-\sqrt{2}f^2_4, \ X_2f^2_8=0,
\ X_2f^2_9=-\sqrt{2}f^2_3,
\end{array}
\right.
\end{eqnarray*}

\begin{eqnarray*}
\left\{
\begin{array}{l}
X_3f^2_1=0, \ X_3f^2_2=-\sqrt{2}f^2_3, \ X_3f^2_3=\sqrt{2}f^2_2, \
X_3f^2_4=0, \ X_3f^2_5=\sqrt{2}f^2_7,
\\ \mbox{} \\
X_3f^2_6=\sqrt{2}f^2_8, \ X_3f^2_7=-\sqrt{2}f^2_5, \
X_3f^2_8=-\sqrt{2}f^2_6, \ X_3f^2_9=0.
\end{array}
\right.
\end{eqnarray*}
From these relations, $J^{\psi}$ clearly preserves the subspaces generated by:
$$
B_1=\{f^2_1d\psi(X_2),f^2_6d\psi(X_2),f^2_8d\psi(X_2)\}\cup
\{f^2_1d\psi(X_3),f^2_6d\psi(X_3),f^2_8d\psi(X_3)\},
$$
$$
B_2=\{-f^2_9d\psi(X_2),f^2_3d\psi(X_2),f^2_2d\psi(X_2)\}\cup
\{-f^2_9d\psi(X_3),f^2_3d\psi(X_3),f^2_2d\psi(X_3)\},
$$
and
$$
B_3=\{f^2_4d\psi(X_2),f^2_7d\psi(X_2),-f^2_5d\psi(X_2)\}\cup
\{f^2_4d\psi(X_3),f^2_7d\psi(X_3),-f^2_5d\psi(X_3)\},
$$
and is represented by the same matrix in any of these three
bases:
$$
\left(
\begin{array}{cccccccccc}
4 & 0 & 0 & 0 & 0 & 4 \\
0 & 4 & 0 & 0 & 0 & 0 \\
0 & 0 & 4 & -4 & 0 & 0 \\
0 & 0 & -4 & 4 & 0 & 0 \\
0 & 0 & 0 & 0 & 4 & 0 \\
4 & 0 & 0 & 0 & 0 & 4 \\
\end{array}
\right),
$$
of eigenvalues $0$, $4$ and $8$, all of multiplicity $2$.

$iv)$ For $k>2$, $J^{\psi}$ restricted to $S^{\psi}_{\lambda_k}$ is positive
definite. Indeed, the equality:
$$
\big(J^{\psi}(f_2d\psi(X_2)+f_3d\psi(X_3)),f_2d\psi(X_2)+f_3d\psi(X_3)\big)=
\int_{\s^3(\sqrt{2})}\Big(\lambda_k\big((f_2)^2+(f_3)^2\big)
-4\sqrt{2}(X_1f_3)f_2\Big) \ v_g,
$$
calls for an upper bound of
$\big\vert \int_{\s^3(\sqrt{2})}(X_1f_3)f_2 \ v_g\big\vert$:

\noindent By Theorem ~\ref{eq:Hilbertbasis}, consider an $L^2$-orthonormal
basis $\{f^k_1,\ldots,f^k_{m_{\lambda_k}}\}$
of $\{f\in C^{\infty}(\s^3(\sqrt{2}))\vert \Delta f=\lambda_kf\}$
such that $\Delta^V f^k_i=c_if^k_i$, for any $i\in
\{1,\ldots,m_{\lambda_k}\}$, $m_{\lambda_k}$ being the
multiplicity of $\lambda_k$. We have:
$$
\Delta^V f_3=\Delta^V\Big(\sum_{i=1}^{m_{\lambda_k}}a^if^k_i\Big)=
\sum_{i=1}^{m_{\lambda_k}}a^ic_if^k_i, \quad a^i\in
\r,
$$
and
\begin{eqnarray*}
\int_{\s^3(\sqrt{2})}(X_1f_3)^2 \ v_g&=&\int_{\s^3(\sqrt{2})}(\Delta^V
f_3)f_3 \ v_g \\
&=&\sum_{i,j=1}^{m_{\lambda_k}}\int_{\s^3(\sqrt{2})}a^ic_ia^jf^k_if^k_j \ v_g
=\sum_{i,j=1}^{m_{\lambda_k}}a^ic_ia^j\delta_{ij}
=\sum_{i=1}^{m_{\lambda_k}}(a^i)^2c_i
\\
&\leq&c\sum_{i=1}^{m_{\lambda_k}}(a^i)^2=c\int_{\s^3(\sqrt{2})}(f_3)^2 \ v_g,
\end{eqnarray*}
where $c=\max\{c_1,\ldots,c_{m_{\lambda_k}}\}$. Plugging in the Cauchy
inequality
$$
\Big\vert \int_{\s^3(\sqrt{2})}(X_1f_3)f_2 \ v_g\Big\vert
\leq
\Big(\int_{\s^3(\sqrt{2})}(X_1f_3)^2 \ v_g\Big)^{\frac{1}{2}}
\Big(\int_{\s^3(\sqrt{2})}(f_2)^2 \ v_g\Big)^{\frac{1}{2}},
$$
yields:
$$
\Big\vert \int_{\s^3(\sqrt{2})}(X_1f_3)f_2 \ v_g\Big\vert
\leq\sqrt{c}
\Big(\int_{\s^3(\sqrt{2})}(f_3)^2 \ v_g\Big)^{\frac{1}{2}}
\Big(\int_{\s^3(\sqrt{2})}(f_2)^2 \ v_g\Big)^{\frac{1}{2}}
\leq\frac{\sqrt{c}}{2}\int_{\s^3(\sqrt{2})}(f_2)^2+(f_3)^2 \ v_g,
$$
and
$$
\big(J^{\psi}(f_2d\psi(X_2)+f_3d\psi(X_3)),f_2d\psi(X_2)+f_3d\psi(X_3)\big)
\geq
(\lambda-2\sqrt{2c})\int_{\s^3(\sqrt{2})}(f_2)^2+(f_3)^2 \ v_g,
$$
hence $J^{\psi}$ is positive definite.

Finally, from the eigenvalues of the above matrices,
we reobtain that $\psi$ has index $4$ and nullity $8$.

\subsection{Description of $\ker J^{\psi}$}

\begin{theorem}
The Hopf map $\psi:\s^3(\sqrt{2})\to \s^2(\sqr)$ is a harmonic map
whose Jacobi operator $J^{\psi}$ is negative definite on
$\{d\psi(\grad f_1)\vert \Delta f_1=\lambda_1 f_1\}\subset
C(\psi^{-1}T\s^2(\sqr))$.
\end{theorem}

\begin{proof}
Recall that $J^{\psi}$ is negative definite on
$$
\{fd\psi(X_2)-\sqrt{2}(X_1f)d\psi(X_3)\vert \Delta f=\lambda_1f\},
$$
so $f(x)=\langle u,x\rangle$, $u\in\r^4$, and $f_1(x)=\langle v,x\rangle$,
where $v=\sqrt{2}(-u^3,u^4,u^1,-u^2)$, satisfies:
$$
fd\psi(X_2)-\sqrt{2}(X_1f)d\psi(X_3)=(X_2f_1)d\psi(X_2)+(X_3f_1)d\psi(X_3)
=d\psi(\grad f_1).
$$
\end{proof}

\begin{proposition}
Let $X=f_2X_2+f_3X_3$ be a horizontal vector field, where
$f_2,f_3\in C^{\infty}(\s^3(\sqrt{2}))$. Then $X$ is
basic if and only if
$$
\Delta^Vf_2=2f_2 \quad \hbox{and} \quad
f_3=-\frac{1}{\sqrt{2}}(X_1f_2).
$$
\end{proposition}

\begin{proof}
As
$d\psi(X_2)=b(f^2_9,f^2_1,-f^2_4)$ and $d\psi(X_3)=b(f^2_2,-f^2_8,-f^2_5)$:
$$
d\psi(f_2X_2+f_3X_3)=
b(f_2f^2_9+f_3f^2_2,f_2f^2_1-f_3f^2_8,-f_2f^2_4-f_3f^2_5),
$$
and $X=f_2X_2+f_3X_3$ is basic if and only if:

$\left\{
\begin{array}{l}
X_1(f_2f^2_9+f_3f^2_2)=0
\\ \mbox{} \\
X_1(f_2f^2_1-f_3f^2_8)=0
\\ \mbox{} \\
X_1(f_2f^2_4+f_3f^2_5)=0.
\end{array}\right.
$

\bigskip

$\Longleftrightarrow\left\{
\begin{array}{l}
(X_1f_2)+\sqrt{2}f_3=0
\\ \mbox{} \\
(X_1f_3)-\sqrt{2}f_2=0
\end{array}\right.
$

and this system is equivalent to $\Delta^Vf_2=2f_2$ and
$f_3=-\frac{1}{\sqrt{2}}(X_1f_2)$.
\end{proof}

This enables a characterization of $\ker J^{\psi}$.

\begin{proposition}
\label{eq:basicfields}
Let $X=f_2X_2+f_3X_3$, $\Delta
f_2=\lambda_kf_2$ and $\Delta f_3=\lambda_k f_3$. Then
$d\psi(X)\in\ker J^{\psi}$ if and only if $k=0$ or $k=2$ and $X$
is basic.
\end{proposition}

\begin{proof}
We saw that, if $d\psi(X)\in\ker J^{\psi}$, $k$ must be $0$ or $2$.

If $k=0$, $f_2$ and $f_3$ are constants, and if $k=2$,
$\ker J^{\psi}_{\vert S^{\psi}_{\lambda_2}}
=\{fd\psi(X_2)-\frac{1}{\sqrt{2}}(X_1f)d\psi(X_3)\vert\Delta^V f=2f\}$.
But $X=fX_2-\frac{1}{\sqrt{2}}(X_1f)X_3$, where $\Delta^Vf=2f$, is
a basic vector field.
\end{proof}

An even more precise description is:

\begin{theorem}
\label{eq:descriptionker}
The kernel of the Jacobi operator $J^{\psi}$ of the Hopf map is:
$$
\{d\psi(X)\vert X\in C(T\s^3(\sqrt{2})), X \ \hbox{is Killing}\}
\oplus \big\{(\grad\tilde{f})\circ\psi \vert \tilde{f}\in
C^{\infty}(\s^2(\sqr)), \Delta\tilde{f}=\mu_1\tilde{f}\big\}.
$$
\end{theorem}

\begin{proof}
First note that:
$$
\{Y\circ\psi \vert Y\in C(T\s^2(\sqr)), Y \ \hbox{is
Killing}\}=\Span\{d\psi(X_4),d\psi(X_5),d\psi(X_6)\}\subset
S^{\psi}_{\lambda_2},
$$
and, as for any $\tilde{f}\in C^{\infty}(\s^2(\sqr))$ we have
$(\grad\tilde{f})\circ\psi=d\psi(\grad(\tilde{f}\circ\psi))$,
$$
\big\{(\grad\tilde{f})\circ\psi \vert \tilde{f}\in
C^{\infty}(\s^2(\sqr)),
\Delta\tilde{f}=\mu_1\tilde{f}\big\}=\big\{d\psi(\grad(\tilde{f}\circ\psi))
\big\}\subset S^{\psi}_{\lambda_2}.
$$
As $X^{\tiny{\hor}}_4$, $X^{\tiny{\hor}}_5$, $X^{\tiny{\hor}}_6$
and $\grad(\tilde{f}\circ\psi)$ are basic,
by Proposition ~\ref{eq:basicfields}:
$$
\{Y\circ\psi \vert Y\in C(T\s^2(\sqr)), Y \ \hbox{is
Killing}\}\subset \ker J^{\psi}
$$
and
$$
\big\{(\grad\tilde{f})\circ\psi \vert \tilde{f}\in
C^{\infty}(\s^2(\sqr)), \Delta\tilde{f}=\mu_1\tilde{f}\big\}
\subset \ker J^{\psi}.
$$
Since both $ \{Y\circ\psi \vert Y\in C(T\s^2(\sqr)), Y \ \hbox{is
Killing}\} $ and $\big\{(\grad\tilde{f})\circ\psi \vert
\tilde{f}\in C^{\infty}(\s^2(\sqr)),
\Delta\tilde{f}=\mu_1\tilde{f}\big\} $ have dimension $3$ and
$\psi$ has nullity $8$:
\begin{eqnarray*}
\ker J^{\psi}&=&\Span\{d\psi(X_2),d\psi(X_3)\} \\
&& \oplus \{Y\circ\psi \vert Y\in C(T\s^2(\sqr)), Y
\ \hbox{is Killing}\} \\
&& \oplus \big\{(\grad\tilde{f})\circ\psi \vert \tilde{f}\in
C^{\infty}(\s^2(\sqr)), \Delta\tilde{f}=\mu_1\tilde{f}\big\}.
\end{eqnarray*}
As
\begin{eqnarray*}
\{d\psi(X)\vert X\in C(T\s^3(\sqrt{2})), X \ \hbox{is Killing}\}&=&
\Span\{d\psi(X_2),d\psi(X_3)\}
\\ &&
\oplus \{Y\circ\psi \vert Y\in C(T\s^2(\sqr)), Y \ \hbox{is
Killing}\},
\end{eqnarray*}
the theorem follows.
\end{proof}

\begin{remark} Theorem ~\ref{eq:descriptionker} can be proven differently:

If $X$ is Killing,
then $\Delta^{\psi}(d\psi(X))=2d\psi(X)$ and so
$J^{\psi}(d\psi(X))=0$.

Furthermore, if $Y\in C(T\s^2(\sqr))$, one computes that:
$$
J^{\psi}(Y\circ\psi)=J^{\bf
1}(Y)\circ\psi=(\Delta_H(Y)-4Y)\circ\psi,
$$
where ${\bf 1}:\s^2(\sqr)\to \s^2(\sqr)$ is the identity map. As
the eigenvalues of $\Delta_H$ on $\s^2(\sqr)$ are
$\{2(k+1)(k+2)\vert k\in \n\}$, $J^{\psi}(Y\circ\psi)=0$ if and
only if $Y$ is Killing or $Y=\grad\tilde{f}$ where $\tilde{f}\in
C^{\infty}(\s^2(\sqr))$ and $\Delta\tilde{f}=\mu_1\tilde{f}$. Now,
a simple dimension count implies Theorem ~\ref{eq:descriptionker}.
\end{remark}

\section{The biharmonic index and nullity of $\psi$}

By definition, any harmonic map is a minimum of the bienergy, thus
biharmonic stable. To ascertain the (biharmonic) nullity of the
Hopf map $\psi:\s^3(\sqrt{2})\to \s^2(\sqr)$, we need:
$$
\trace\langle V,d\psi\cdot\rangle d\psi\cdot=V, \quad \vert
d\psi\vert^2=2,
$$
$$
\trace\langle d\psi\cdot,\Delta^{\psi}V\rangle
d\psi\cdot=\Delta^{\psi}V, \quad \trace\langle
d\psi\cdot,\trace\langle d\psi\cdot,V\rangle
d\psi\cdot\rangle d\psi\cdot=V,
$$
and, replacing in ~\eqref{eq:J1}, we obtain:
\begin{eqnarray*}
I^{\psi}(V)&=&\Delta^{\psi}\Delta^{\psi}(V)-4\Delta^{\psi}V+4V \\
&=&J^{\psi}(J^{\psi}(V)).
\end{eqnarray*}

Again, since homothetic transformations do not alter the spectrum of the
operator $I^{\psi}$, we attain:

\begin{theorem}
The Hopf map $\psi:\s^3\to \s^2$
is biharmonic stable and its biharmonic and harmonic nullities coincide,
i.e. equal $8$.
\end{theorem}

\section{On the biharmonic index and nullity of $\phi$}

Recall that, if $\psi:M\to \s^n(\sqr)$ is a nonconstant harmonic
map of constant energy density, then $\phi={\bf i}\circ\psi:M\to
\s^{n+1}$ is nonharmonic biharmonic. As to its stability we have:

\begin{theorem}
Let $M$ be a compact manifold and $\psi:(M,g)\to \s^n(\sqr)$ a
nonconstant harmonic map of constant energy density. The map
$\phi={\bf i}\circ \psi:M\to \s^{n+1}$ is biharmonic unstable.
\end{theorem}

\begin{proof}
Simple considerations furnish:
$$
\trace\langle\eta,d\phi\cdot\rangle d\phi\cdot=0, \quad
\vert d\phi\vert^2=2e(\psi), \quad
\langle d\tau(\phi),d\phi\rangle\eta=-4(e(\psi))^2\eta,
$$
$$
\vert\tau(\phi)\vert^2=4(e(\psi))^2, \quad
\trace\langle\eta,d\tau(\phi)\cdot\rangle d\phi\cdot=0, \quad
\trace\langle\tau(\phi),d\eta\cdot\rangle d\phi\cdot=0,
$$
$$
\langle\tau(\phi),\eta\rangle\tau(\phi)=4(e(\psi))^2\eta, \quad
\trace\langle d\phi\cdot,\Delta^{\phi}\eta\rangle d\phi\cdot=
(\Delta^{\phi}\eta)^{\top}, \quad
\langle d\eta,d\phi\rangle\tau(\phi)=-4(e(\psi))^2\eta,
$$
and, by ~\eqref{eq:J1}:
$$
(I^{\phi}(\eta),\eta)=\int_M\big(\vert\Delta^{\phi}\eta\vert^2
-4e(\psi)\langle\Delta^{\phi}\eta,\eta\rangle-12(e(\psi))^2\big) \ v_g.
$$
Since $\Delta^{\phi}\eta=2e(\psi)\eta$,
$(I^{\phi}(\eta),\eta)$ is strictly negative, i.e. $\phi$ is
biharmonic unstable.
\end{proof}

For the Hopf map, we set $V=d\phi(X)$, where $X\in C(T\s^3(\sqrt{2}))$ and,
to find $I^{\phi}(V)$, evaluate the terms:
$$
\trace\langle V,d\phi\cdot\rangle d\phi\cdot=V, \quad \vert
d\phi\vert^2=2, \quad \langle d\tau(\phi),d\phi\rangle V=-4V, \quad
\vert\tau(\phi)\vert^2V=4V,
$$
$$
\trace\langle V,d\tau(\phi)\cdot\rangle d\phi\cdot=-2V, \quad
\trace\langle\tau(\phi),dV\cdot\rangle d\phi\cdot=2V, \quad
\langle\tau(\phi),V\rangle\tau(\phi)=0,
$$
$$
\trace\langle d\phi\cdot,\Delta^{\phi}V\rangle
d\phi\cdot=(\Delta^{\phi}V)^{\top}, \quad \langle
dV,d\phi\rangle\tau(\phi)=\trace\langle\nabla^{\psi}d\psi(X),
d\psi\rangle\tau(\phi).
$$
Replacing in ~\eqref{eq:J1}, we reach:
$$
I^{\phi}(V)=\Delta^{\phi}\Delta^{\phi}V-3\Delta^{\phi}V
+(\Delta^{\phi}V)^{\top}-3V
+2\trace\langle\nabla^{\psi}d\psi(X),d\psi\rangle\tau(\phi).
$$
Now, $\Delta^{\phi}V$ can be written
$$
\Delta^{\phi}V=\Delta^{\psi}V+2\trace\langle\nabla^{\psi}d\psi(X),
d\psi\rangle\eta+V,
$$
thus, if $f_2,f_3\in C^{\infty}(\s^2(\sqrt{2}))$ with
$\Delta f_2=\lambda_kf_2$ and $\Delta f_3=\lambda_kf_3$,
$$
\Delta^{\psi}(d\psi(f_2X_2))=(\lambda_k+2)f_2d\phi(X_2)
+2\sqrt{2}(X_1f_2)d\phi(X_3),
$$
$$
\Delta^{\psi}(d\psi(f_3X_3))=(\lambda_k+2)f_3d\phi(X_3)
-2\sqrt{2}(X_1f_3)d\phi(X_2),
$$
$$
\trace\langle\nabla^{\psi}d\psi(f_2X_2),d\psi\rangle\eta=(X_2f_2)\eta,
\
\trace\langle\nabla^{\psi}d\psi(f_3X_3),d\psi\rangle\eta=(X_3f_3)\eta,
$$
therefore
\begin{eqnarray*}
I^{\phi}(f_2d\phi(X_2))&=&\Big(\big( (\lambda_k)^2+4\lambda_k\big
)f_2-8(X_1X_1f_2)-4(X_2X_2f_2)\Big)d\phi(X_2) \\
&&+\big(4\sqrt{2}(\lambda_k+2)(X_1f_2)-4(X_3X_2f_2)\big)d\phi(X_3)
\\
&&+\big(4\lambda_k(X_2f_2)+4\sqrt{2}(X_3X_1f_2)\big)\eta,
\end{eqnarray*}
\begin{eqnarray*}
I^{\phi}(f_3d\phi(X_3))&=&\big(-4\sqrt{2}(\lambda_k+2)(X_1f_3)
-4(X_2X_3f_3)\big)d\phi(X_2)
\\
&&+\Big(\big((\lambda_k)^2+4\lambda_k\big)f_3-8(X_1X_1f_3)-4(X_3X_3f_3)
\Big)d\phi(X_3) \\
&&+\big(4\lambda_k(X_3f_3)-4\sqrt{2}(X_2X_1f_3)\big)\eta.
\end{eqnarray*}

Set $W=f\eta$, where $f\in C^{\infty}(\s^3(\sqrt{2}))$, accordingly:
$$
I^{\phi}(W)=\Delta^{\phi}\Delta^{\phi}W-4\Delta^{\phi}W
+(\Delta^{\phi}W)^{\top}-12W+4d\phi(\grad f)
$$
and
$$
\Delta^{\phi}W=(\Delta f+2f)\eta-2(X_2f)d\phi(X_2)-2(X_3f)d\phi(X_3).
$$

If $f\in C^{\infty}(\s^2(\sqrt{2}))$ with
$\Delta f=\lambda_kf$:
\begin{eqnarray*}
I^{\phi}(f\eta)&=&\big(-4\lambda_k(X_2f)+4\sqrt{2}(X_1X_3f)
\big)d\phi(X_2) \\
&&+\big(-4\lambda_k(X_3f)-4\sqrt{2}(X_1X_2f)\big)d\phi(X_3) \\
&&+\Big(\big((\lambda_k)^2-16\big)f-4(X_2X_2f)-4(X_3X_3f)\Big)\eta.
\end{eqnarray*}
As $X_1$, $X_2$ and $X_3$ preserve the eigenspaces of the Laplacian,
$I^{\phi}$ preserves the subspace
$$
S^{\phi}_{\lambda_k}=\{f_2d\phi(X_2)\vert \Delta
f_2=\lambda_kf_2\}\oplus \{f_3d\phi(X_3)\vert \Delta
f_3=\lambda_kf_3\}\oplus \{f\eta\vert \Delta f=\lambda_kf\},
$$
for any $k\in \n$.

In order to determine the index of $I^{\phi}$ we follow the same approach as
for $J^{\psi}$.

$i) \ k=0.$ In this case
$$
B^{\phi}_{\lambda_0}=\big\{\frac{1}{c}d\phi(X_2),\frac{1}{c}d\phi(X_3),
\frac{1}{c}\eta\big\}
$$
is an $L^2$-orthonormal basis for $S^{\phi}_{\lambda_0}$,
where $c^2=\vol \big(\s^3(\sqrt{2})\big)$. As
$I^{\phi}(d\phi(X_2))=I^{\phi}(d\phi(X_3))=0$ and
$I^{\phi}(\eta)=-16\eta$, the matrix of $I^{\phi}$ in
$B^{\phi}_{\lambda_0}$ is
$$
\left(
\begin{array}{ccc}
0 & 0 & 0 \\
0 & 0 & 0 \\
0 & 0 & -16
\end{array}
\right).
$$

$ii) \ k=1.$ An $L^2$-orthonormal basis for $S^{\phi}_{\lambda_1}$ is
$$
B^{\phi}_{\lambda_1}=\{f^1_id\phi(X_2)\}_{i=1}^4\cup
\{f^1_jd\phi(X_3)\}_{j=1}^4\cup \{f^1_k\eta\}_{k=1}^4.
$$
Computing
$I^{\phi}(f^1_id\phi(X_2))$, $I^{\phi}(f^1_jd\phi(X_2))$ and
$I^{\phi}(f^1_kd\phi(X_2))$, $i,j,k\in\{1,\ldots, 4\}$ supplies the
matrix
$$
\left(
\begin{array}{cccccccccccc}
\frac{57}{4} & 0 & 0 & 0 & 0 & -12 & 0 & 0 & 0 & 0 & -\sqrt{2} & 0
\\
0 & \frac{57}{4} & 0 & 0 & 12 & 0 & 0 & 0 & 0 & 0 & 0 & \sqrt{2}
\\
0 & 0 & \frac{57}{4} & 0 & 0 & 0 & 0 & -12 & \sqrt{2} & 0 & 0 & 0
\\
0 & 0 & 0 & \frac{57}{4} & 0 & 0 & 12 & 0 & 0 & -\sqrt{2} & 0 & 0
\\
0 & 12 & 0 & 0 & \frac{57}{4} & 0 & 0 & 0 & 0 & 0 & 0 & -\sqrt{2}
\\
-12 & 0 & 0 & 0 & 0 & \frac{57}{4} & 0 & 0 & 0 & 0 & -\sqrt{2} & 0
\\
0 & 0 & 0 & 12 & 0 & 0 & \frac{57}{4} & 0 & 0 & \sqrt{2} & 0 & 0
\\
0 & 0 & -12 & 0 & 0 & 0 & 0 & \frac{57}{4} & \sqrt{2} & 0 & 0 & 0
\\
0 & 0 & \sqrt{2} & 0 & 0 & 0 & 0 & \sqrt{2} & -\frac{39}{4} & 0 &
0 & 0 \\
0 & 0 & 0 & -\sqrt{2} & 0 & 0 & \sqrt{2} & 0 & 0 & -\frac{39}{4} &
0 & 0 \\
-\sqrt{2} & 0 & 0 & 0 & 0 & -\sqrt{2} & 0 & 0 & 0 & 0 &
-\frac{39}{4} & 0 \\
0 & \sqrt{2} & 0 & 0 & -\sqrt{2} & 0 & 0 & 0 & 0 & 0 & 0 &
-\frac{39}{4}
\end{array}
\right)
$$
of eigenvalues $-\frac{15}{4}+2\sqrt{10}$, $-\frac{15}{4}-2\sqrt{10}$
and $\frac{105}{4}$, all three of multiplicity $4$.

$iii) \ k=2.$
$
B^{\phi}_{\lambda_2}=\{f^2_id\phi(X_2)\}_{i=1}^9\cup
\{f^2_jd\phi(X_3)\}_{j=1}^9\cup \{f^2_k\eta\}_{k=1}^9
$
is an $L^2$-orthonormal basis of $S^{\phi}_{\lambda_2}$. As in the
harmonic case, it is enough to compute the matrix of
$I^{\phi}$ in
$$
\{f^2_id\phi(X_2)\}_{i\in\{1,6,8\}}\cup
\{f^2_jd\phi(X_3)\}_{j\in\{1,6,8\}}\cup
\{f^2_k\eta\}_{k\in\{1,6,8\}}.
$$
This matrix is
$$
\left(
\begin{array}{ccccccccc}
56 & 0 & 0 & 0 & 0 & 40 & 0 & 8\sqrt{2} & 0 \\
0 & 40 & 0 & 0 & 0 & 0 & -16\sqrt{2} & 0 & 0 \\
0 & 0 & 48 & -48 & 0 & 0 & 0 & 0 & 0 \\
0 & 0 & -48 & 48 & 0 & 0 & 0 & 0 & 0 \\
0 & 0 & 0 & 0 & 40 & 0 & 0 & 0 & 16\sqrt{2} \\
40 & 0 & 0 & 0 & 0 & 56 & 0 & -8\sqrt{2} & 0 \\
0 & -16\sqrt{2} & 0 & 0 & 0 & 0 & 8 & 0 & 0 \\
8\sqrt{2} & 0 & 0 & 0 & 0 & -8\sqrt{2} & 0 & 16 & 0 \\
0 & 0 & 0 & 0 & 16\sqrt{2} & 0 & 0 & 0 & 8
\end{array}
\right)
$$
and its eigenvalues are $0$, $96$, $24+16\sqrt{3}$, $24-16\sqrt{3}$, of
multiplicity two, and $32$ of multiplicity one.

$iv) \ k>2.$ Let $V=f_2d\phi(X_2)+f_3d\phi(X_3)+f\eta$, where
$\Delta f_2=\lambda_kf_2$, $\Delta f_3=\lambda_kf_3$, $\Delta
f=\lambda_kf$. By a straightforward computation:
\begin{eqnarray*}
(I^{\phi}(V),V)&=& \int_{\s^3(\sqrt{2})} \Big(\big(
(\lambda_k)^2+4\lambda_k\big)(f_2)^2+
\big((\lambda_k)^2+4\lambda_k\big)(f_3)^2+\big((\lambda_k)^2-16\big)f^2
\\
&&
\ +4(X_2f_2+X_3f_3)^2+\big(2\sqrt{2}(X_1f_2)-2(X_3f)\big)^2
+\big(2\sqrt{2}(X_1f_3)+2(X_2f)\big)^2
\\
&& \
-8\sqrt{2}(\lambda_k+2)(X_1f_3)f_2-8\lambda_k(X_2f)f_2
-8\lambda_k(X_3f)f_3\Big) \ v_g
\\
&\geq&
\int_{\s^3(\sqrt{2})} \Big(\big(
(\lambda_k)^2+4\lambda_k\big)(f_2)^2+
\big((\lambda_k)^2+4\lambda_k\big)(f_3)^2+\big((\lambda_k)^2-16\big)f^2
\\
&& \
-8\sqrt{2}(\lambda_k+2)(X_1f_3)f_2-8\lambda_k(X_2f)f_2
-8\lambda_k(X_3f)f_3\Big) \ v_g
\end{eqnarray*}
As
$$
\Big\vert
\int_{\s^3(\sqrt{2})}(X_if_2)f_3 \ v_g\Big\vert\leq\frac{\sqrt{c}}{2}
\int_{\s^3(\sqrt{2})}(f_2)^2+(f_3)^2 \ v_g,
$$
we obtain
\begin{eqnarray*}
(I^{\phi}(V),V)&\geq&\int_{\s^3(\sqrt{2})}\Big(\big(
(\lambda_k)^2+4\lambda_k-4\sqrt{2c}(\lambda_k+2)
-4\lambda_k\sqrt{c}\big)\big((f_2)^2+(f_3)^2\big)
\\
&& \ +\big((\lambda_k)^2-8\lambda_k\sqrt{c}-16\big)f^2\Big) \ v_g.
\end{eqnarray*}
Let
$
A=(\lambda_k)^2+4\lambda_k-4\sqrt{2c}(\lambda_k+2)-4\lambda_k\sqrt{c}
$
and
$
B=(\lambda_k)^2-8\lambda_k\sqrt{c}-16,
$
then
$$
A\geq\frac{k^2(k+2)^2}{4}+2k(k+2)-2k\big(k(k+2)+4\big)-\sqrt{2}k^2(k+2)
>0, \quad \forall k\geq 12
$$
and
$$
B\geq\frac{k^2(k+2)^2}{4}-2\sqrt{2}k^2(k+2)-16>0, \quad \forall k\geq
10.
$$
Therefore, $I^{\phi}$ restricted to
$S^{\phi}_{\lambda_k}$ is positive definite for any $k\geq 12$.

From this analysis we conclude:

\begin{theorem}
The index of the biharmonic map $\phi:\s^3\to \s^3$ is at least $11$
while its nullity is bounded from bellow by $8$.
\end{theorem}

\subsection{The basic biharmonic index and nullity of $\phi$.}
Let $S$ be the subspace of $C(\phi^{-1}T\s^3)$ defined by:
$$
S=\{X\circ\psi\vert X\in C(T\s^2(\sqr))\}\oplus
\{(\tilde{f}\circ\psi)\eta\vert \tilde{f}\in
C^{\infty}(\s^2(\sqr))\}.
$$
Set $V=X\circ\psi$ and $W=f\eta$, where $f=\tilde{f}\circ\psi$. We have:
$$
I^{\phi}(V)=\Delta^{\phi}\Delta^{\phi}V-3\Delta^{\phi}V+(\Delta^{\phi}
V)^{\top}-3V-4((\Div X)\circ\psi)\eta,
$$
$$
I^{\phi}(W)=\Delta^{\phi}\Delta^{\phi}W-4\Delta^{\phi}W
+(\Delta^{\phi}W)^{\top}-12W+4(\grad\tilde{f})\circ\psi,
$$
$$
\Delta^{\phi}V=(X-\trace\nabla^2X)\circ\psi+2((\Div
X)\circ\psi)\eta=(\Delta^{{\bf i}}X)\circ\psi,
$$
$$
\Delta^{\phi}(f\eta)=((\Delta\tilde{f}+2\tilde{f})\circ\psi)\eta-2(\grad
\tilde{f})\circ\psi=(\Delta^{{\bf i}}(\tilde{f}\eta))\circ\psi,
$$
where ${\bf i}:\s^2(\sqr)\to \s^3$ is the canonical inclusion.
From these relations we deduce that $ I^{\phi}(V)=(I^{\bf
i}(X))\circ\psi $ and $ I^{\phi}(W)=(I^{\bf
i}(\tilde{f}\eta))\circ\psi, $ so $I^{\phi}$ preserves the
subspace $S$ and the restriction of $I^{\phi}$ to $S$,
$I^{\phi}_{\vert S}$, has the same spectrum as $I^{\bf i}$. Thus,
using Propositions $5.1.$ and $5.2.$ in ~\cite{ELCO} we conclude:

\begin{theorem}
The index of $I^{\phi}_{\vert S}$ is $1$ and its nullity is $6$.
Moreover, the subspace which gives
the index of $I^{\phi}_{\vert S}$ is $\{c\eta\vert c\in\r\}$ and the
subspace which gives
the nullity is $\{X\circ\psi\vert X \ \hbox{is Killing}\}\oplus
\{2(\tilde{f}\circ\psi)\eta+(\grad \tilde{f})\circ\psi\vert
\Delta\tilde{f}=\mu_1\tilde{f}\}$.
\end{theorem}

\end{document}